\documentclass[fleqn,10pt]{wlscirep}
\usepackage[utf8]{inputenc}
\usepackage[T1]{fontenc}
\title{Assembly along lines in boundary-driven dynamical system}

\author[1,*]{Kulveer Singh}
\author[1,+]{Yitzhak Rabin}
\affil[1]{Department of Physics, and Institute of Nanotechnology and Advanced Materials, Bar-Ilan University,
Ramat Gan 52900, Israel}

\affil[*]{kulveersingh85@gmail.com}
\affil[+]{rabin@biu.ac.il}

\begin{abstract}
We introduce a simple dynamical rule in which each particle locates a particle that is farthest from it and moves towards it. Repeated application of this algorithm results in the formation of unusual dynamical patterns: during the process of assembly the system self-organizes into slices of low particle density separated by lines of increasingly high particle density along which most particles move. As the process proceeds, pairs of lines meet and merge with each other until a single line remains and particles move along it towards the zone of assembly. We show that this pattern is governed by particles (attractors) situated on the instantaneous outer boundary of the system and that both in two and in three dimensions the lines are formed by zigzag motion of a particle towards a pair of nearly equidistant attractors. This novel line-dominated assembly is very different from the local assembly in which particles that move towards their nearest neighbors produce point-like clusters that coalesce into new point-like clusters, etc. 
\end{abstract}
\begin{document}

\flushbottom
\maketitle
%
%
\thispagestyle{empty}


\section*{Introduction}

Systems composed of a large number of autonomous agents interacting with simple rules often exhibit emergent large scale behavior\cite{vicsek}. In living systems such as a bacterial swarms\cite{Czir, Zhang, Cisneros}, flocks of birds\cite{Ballerini, Emlen}, schools of fishes\cite{parrish, Becco}, etc.\cite{couzin, Hayakawa}, a highly coordinated movement among the agents arises in the absence of centralized control, due to the action of individual agents based on the perception of their local environment (i.e the behavior of each agent is determined by that of its neighbors). These synchronized motions observed in natural systems were successfully modeled by various research groups using algorithms based on local behavioral rules of autonomous agents in computer simulations\cite{kunz, vicsek1995, gueron, Zumaya}. For example, using three simple rules, viz. collision avoidance, velocity matching, and flock centering in his `BOID' model\cite{Reynolds}, Craig Reynolds has simulated a coordinated group movement in a flock of birds.

In this paper we propose and analyze a simple non-local algorithm for aggregation of agents according to which, at each moment, every agent senses the locations of all the other agents and moves towards of the farthest agent from it.  This rule implies that the behavior of each agent in a swarm is determined only by agents located at the outer boundary of the swarm. Continued application of this simple rule results in the appearance of anisotropic dynamical patterns composed of low density ``slices'' separated by high density lines. As the system contracts, particles in the slices migrate towards the lines that separate between neighboring slices and continue to move along them towards a gathering point whose position is close to but not coincident with the center of mass of the initial system. In the course of contraction, the number of the slices and of the associated lines decreases due to their coalescence. We show that the dynamic patterns produced by this simple non-local algorithm are qualitatively different from those that arise using a local rule, where an agent moves towards its nearest neighbor.

\section*{Model}
We consider a swarm of $N$ autonomous agents initially randomly distributed in a region. Each agent is modeled as a point particle whose position is updated according to the following simple rule: 
For each particle ($i$) find the particle ($j$) that is furthest  from it at this time and move particle $i$ by distance $\Delta x$ towards particle $j$. 
If there are several particles whose distance from particle $i$  is the same, randomly choose one of them and move by distance $\Delta x$ towards it (while possible in principle, such exact degeneracy was never observed in our simulations).

The above steps are performed by every particle in each iteration (time-step) and therefore, the choice of the farthest particle may change with time. According to this algorithm, the interaction between the particles is not always reciprocal, in the sense that if  particle $j$ is farthest from particle $i$, it is not necessarily the case that $i$ is farthest from $j$. Also, since each particle is affected only by the one particle that is farthest from it and is therefore located at the outer periphery of the system, every particle will move towards the far boundary. The combined effect of such displacements of all the particles (including the boundary particles themselves) towards the far boundary, results in the contraction of the system and the assembly of all the particles in it. The time $t$ is equal to number of discrete time-steps starting from $t=0$.

\section*{Results and Discussion}

\begin{figure*}[ht]
\includegraphics[width=\linewidth]{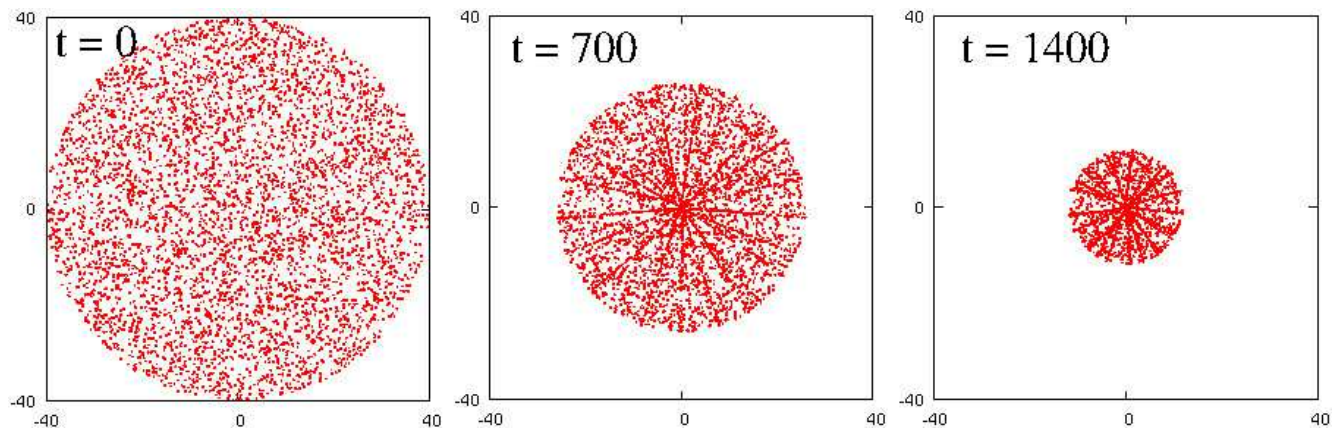}
	\caption{\label{fig:collapse}  Snapshots of a system of 5000 particles at three different times. At $t=0$, particles are randomly distributed inside a circular region. Particles align themselves to form lines in the system($t=700$) and density of particles on the line increases as time progresses($t=1400$). }
\end{figure*}

\begin{figure}[h]
\includegraphics[width=\linewidth]{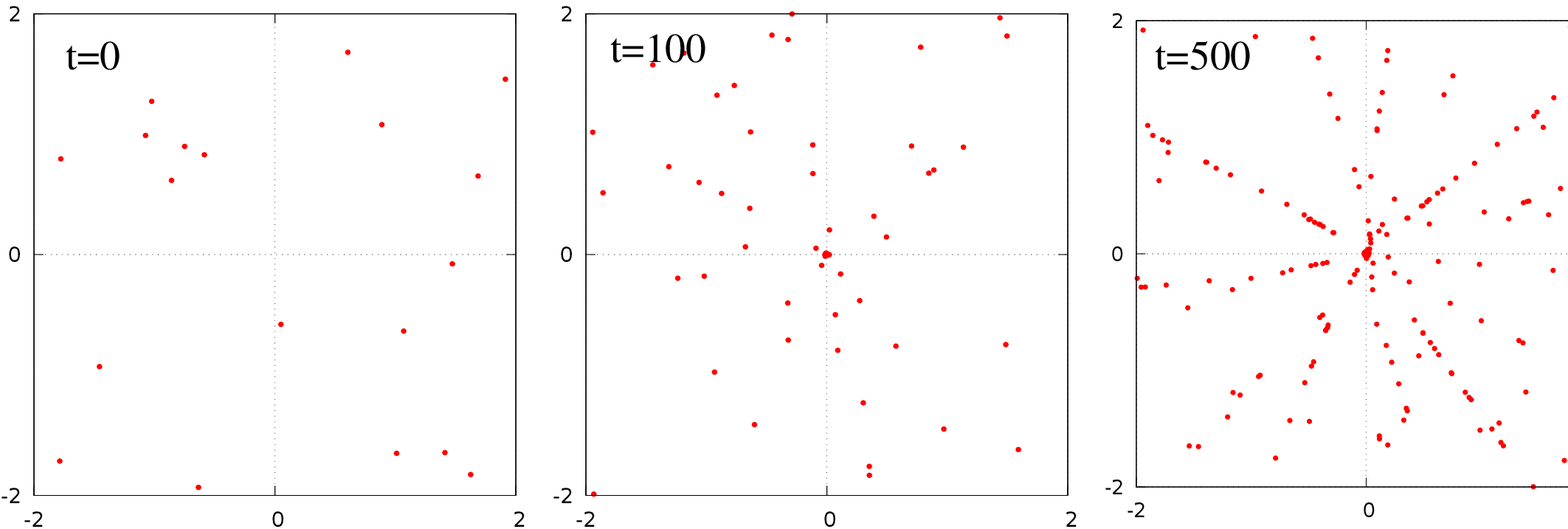}
	\caption{\label{fig:collapse_zoom} Snapshots of the system near the collapse zone at six different times (expanded view). The position of the collapse zone is determined at an early stage of evolution ($t=100$). Several lines meet at the collapse zone; the other lines appear to branch out of these lines ($t=1500$). Density of particles on the lines increases as time progresses (see $t=500, 1500$). Towards the end of the collapse, the number of lines decreases to three ($t=2000$) and finally all the particles assemble in the collapse zone ($t=2500$).    }
\end{figure}

We began the simulation by randomly placing $N=5000$ particles inside a circular disc of radius $R=\sqrt{N/(\pi\rho)}$ with uniform density $\rho=1/\sigma^2$ and chose $\sigma$ as the unit to measure distances($\sigma =1$), and the displacement step $\Delta x=0.02$. Fig. \ref{fig:collapse} shows snapshots of the system at three different times (also see movie M1 in SI). As evident from the snapshots, while initially (at $t=0$) the distribution of particles in the disc is uniform and isotropic, the distribution becomes anisotropic as the system evolves and the rotational symmetry is spontaneously broken. Thus, as particles move inward and the system contracts, it self-organizes into slices of low particle density separated by lines of high density of particles. The formation of lines begins quite early and the density of particles within the lines increases with time. Note that the motion of neighbouring particles along each line is strongly correlated despite the fact that our algorithm does not allow the particles to sense their local environment.   The radius of the circular disk shrinks and all the particles move towards the central region of the disk, as time progresses. Eventually, all the particles assemble in the assembly zone near the center of the disk, defined as a region of width $2\Delta x$ in which all the particles are assembled at the end of the process. 
Fig. \ref{fig:collapse_zoom} shows snapshots of a small section of the system around the assembly zone (see movie M2 in SI). As the system evolves, particles begin assembling in the assembly zone which remains almost fixed till the end of the process. We also observe that only very few lines merge directly at the assembly zone, with the other lines branching out of these lines. 
The total number of lines in the system decreases as the circular disk shrinks. Towards the end of the assembly process the number of lines decreases to three and then to one, and finally all the particles accumulate in the assembly zone (see  Fig. \ref{fig:collapse_zoom}). Furthermore, every point enters the assembly zone along these lines only (see movie M2 in SI). The points in the low density slices between the lines, join the lines as they approach the assembly zone. In order to make sure that the above picture of the dynamics is robust, we repeated the  simulations with many different initial conditions and did not find any qualitative differences in the assembly process.  

To confirm that the formation of lines  is not limited to circular geometry, we performed simulations of random particle distributions in other initial geometries i.e., square and semi-circular domains, and observed four lines in square geometry and only one line (initially) in semi-circular geometry in all the simulations (see movies M3 and M4 in SI). Contrary to these systems, the initial number of lines in the circular disc geometry was observed to depend, albeit weakly, on the initial conditions. This makes the circular disc geometry case more interesting and complex as compared to other geometries.  We also checked the dynamics for larger displacements $\Delta x= 2$ and $5$ (i.e., larger than the mean initial interparticle distance $\sigma=1$ ) and found that the number of lines did not depend on the choice of the stepsize but that the width of these lines increased with increasing $\Delta x$ . In the remainder of the paper, we systematically explore the mechanisms behind the formation of the dynamical patterns observed in the circular disc geometry.. 

\begin{figure}[h]
\includegraphics[width=0.8\linewidth]{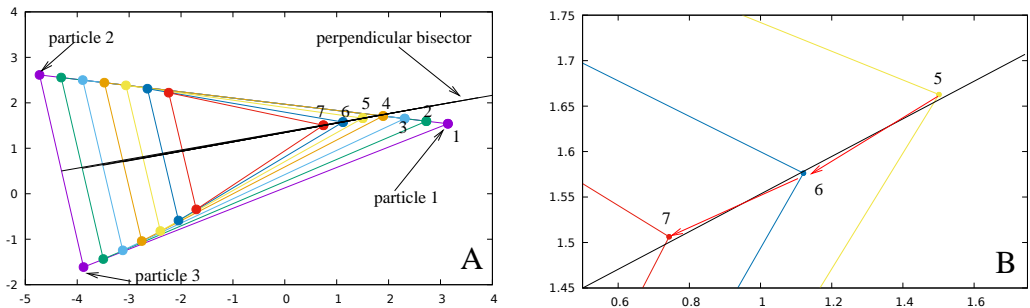}
	\caption{\label{fig:3particles} (A) Three particle system:  differently colored vertices of the triangles represent the positions of the particles at different instants of time. The positions of particle 1 at various times are labeled as $1,2,....,7$. (B) Expanded view of the trajectory of particle 1 which executes a zigzag motion along the perpendicular bisector of particle 2 and particle 3 as it comes within $\Delta x$ from this line (see positions $5,6,7$ of particle 1).  }
\end{figure}

In order to understand the formation of the lines, we turn to examine the dynamics of a smaller system.  Since the two particle case is trivial as  both particles simply move towards each other along the straight line joining them, we consider a three particle system. 
At any instant the positions of the three particles can be thought of as the vertices of a triangle, and in our algorithm the lengths of the sides of this triangle determine the direction of motion of the particles. The two particles which form the longest side of the triangle move towards each other and the third particle moves along the second longest side. 
In Fig. \ref{fig:3particles}A, we plot the positions of the three particles at different times (different colors represent different time instants). We observe that starting from any triangle, the system reaches a stage where three particles form a quasi-isosceles triangle in which two of the longer sides have nearly equal lengths (the lengths differ by less than the step size $\Delta x$). Since at time 1, particle 1 is the farthest from both particles 2 and 3, it acts as an attractor for these particles and they move towards it. Turning our attention to the motion of particle 1, we observe at time instants 1, 2 and 3 in Fig. \ref{fig:3particles}A, particle 1 moves towards particle 2 until it comes within $\Delta x$ of the perpendicular bisector and forms a quasi-isosceles triangle with particles 2 and 3 (also see movie M5 in SI).  Once particle 1 enters this region, the difference between the distances of other two particles from it becomes of the order of step size ($\Delta x$) at which point it has two nearly equidistant attractors. In the next one or two steps, particle 1 crosses the perpendicular bisector and the farthest particle from it becomes particle 3. Now particle 1 moves towards particle 3 and again crosses the perpendicular bisector at which time the farthest particle from it again becomes particle 2. This frequent switching between the two attractors continues and leads to zigzag motion of particle 1 about the perpendicular bisector (see Fig. \ref{fig:3particles}B). Since the amplitude of the zigzag motion is of the order of step size which we have chosen to be very small ($\Delta x =0.02$) compared to the average interparticle distance, particle 1 appears to move along a straight line which is the perpendicular bisector of the side joining the two nearly equidistant attractors, particles 2 and 3.  

Having understood that zigzag motion of a particle along the bisector of two neighboring nearly equidistant attractors in a three particle system appears (at sufficiently low resolution) as motion along a line, we return to the large system case (5000 particles) in the circular disc geometry. A careful examination of the dynamics of assembly shows that while particles on the lines execute zigzag motion (changing direction abruptly and very frequently), all particles inside a ``slice'' located between neighboring  lines move smoothly towards a common point of convergence (see movie M6 in SI where velocity vectors of particles evolving with time are shown). 
This concurs with the expectation that while particles in the interior of a slice move towards a common attractor located near the far outer boundary of the system, those on the boundary line between two neighboring slices execute a zigzag motion whose direction alternates between one of the two nearly equidistant attractors and therefore oscillates around the bisector to the imaginary line connecting these attractors.  The fact that many particles move along the same line indicates that these particles have a common nearly equidistant pair of attractors. The number of lines is identical to the number of  nearly equidistant pairs of attractors in the system. The particles in the slices between the lines move closer to the lines during the process of contraction and eventually join these lines before entering the assembly zone.    

\begin{figure}[h]
\includegraphics[width=\linewidth]{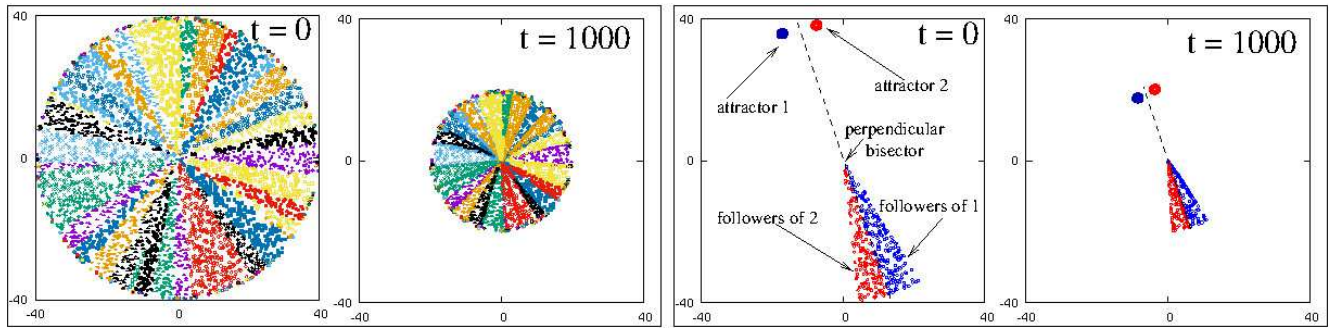}
	\caption{\label{fig:collapse1} Followers of each attractor (towards which they move), are represented by points of the same color. With this convention the system splits into the slices of different colors where each slice corresponds to a different attractor. In the right panel we show two neighboring slices and their corresponding attractors. The boundary between two slices defines a line that lies on the perpendicular bisector of two neighboring attractors. }
\end{figure}

According to our dynamical rules, at any instant of time each particle in the system moves towards another particle (its attractor) and therefore can be termed as the follower of this attractor. Note that all the attractors are located in a narrow annular region close to the outer boundary of the system at time $t$ and therefore the total number of attractors $N_A(t)$ is much smaller that the number of particles in the system $N$. Since each attractor is a follower of some other attractor, the total number of followers is $N$. In order to visualize the way in which the system separates into groups of followers of different attractors, the followers of different attractors are shown in different color in Fig. \ref{fig:collapse1}. This plotting scheme divides the system into differently colored slices where each slice is defined by the followers of one attractor (see movie M7 in SI). Thus, the number of such slices is equal to the number of attractors and the area of each slice is proportional to the number of followers of the corresponding attractor. The border line between neighbouring slices is formed by the perpendicular bisector of the line joining their attractors (see lower panel in Fig. \ref{fig:collapse1} where only two neighbouring slices are shown).

\begin{figure}[h]
\includegraphics[width=0.5\linewidth]{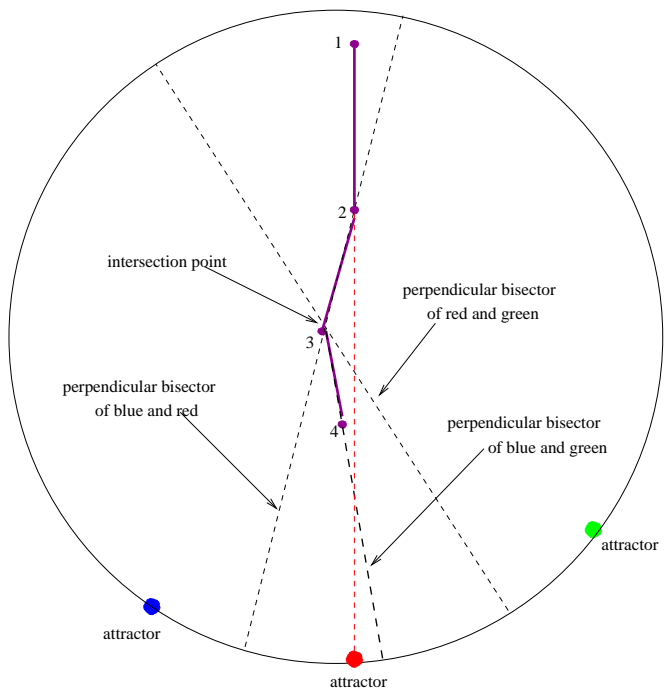}
	\caption{\label{fig:attractor_disappearance} Schematic view of the dynamics of a typical follower of the middle (red) attractor.  }
\end{figure}
Inspection of movie M7 in the SI shows that the number of attractors and that of the corresponding slices of followers and of the associated border lines, decreases with time. The way it happens is demonstrated  in movie M8 in SI where we show three neighboring blue, red and green attractors and their corresponding slices of followers and observe that the area of the middle slice decreases with time and eventually disappears completely, together with its attractor. In order to understand the mechanism behind the decrease of $N_A(t)$ with time, in  Fig. \ref{fig:attractor_disappearance} we present a simplified picture of the above process in which we freeze the attractors and schematically describe the motion of a particle in the slice that corresponds to the middle (red) attractor. As shown in the figure, as the particle moves from point 1 towards its attractor, it reaches the boundary line which is the perpendicular bisector of the imaginary line between the blue and the red attractors (point 2). From this time on it continues to move (in a zigzag fashion) along this line until it approaches the  intersection point (point 3) of two neighbouring perpendicular bisectors which is the circumcenter of the triangle formed by the corresponding three attractors and is therefore equidistant from the three of them. This threefold degeneracy is removed as the particle moves  away from the  intersection point by a small amount $\sim \Delta x$ and the distance between the particle and the red attractor becomes smaller that the distance to the blue and the green attractors. From this time on (point 4), the red point on the boundary stops being an attractor for the particle and the  blue and green points become its quasi-degenerate pair of attractors. The two lines formed by the two neighboring pairs merge to form another line which is the perpendicular bisector of the imaginary line joining the blue and green attractors. Similar dynamics takes place for all the particles in the central slice in Fig. \ref{fig:attractor_disappearance}, leading eventually to the disappearance of the red attractor.

\begin{figure}[h]
\includegraphics[width=0.6\linewidth]{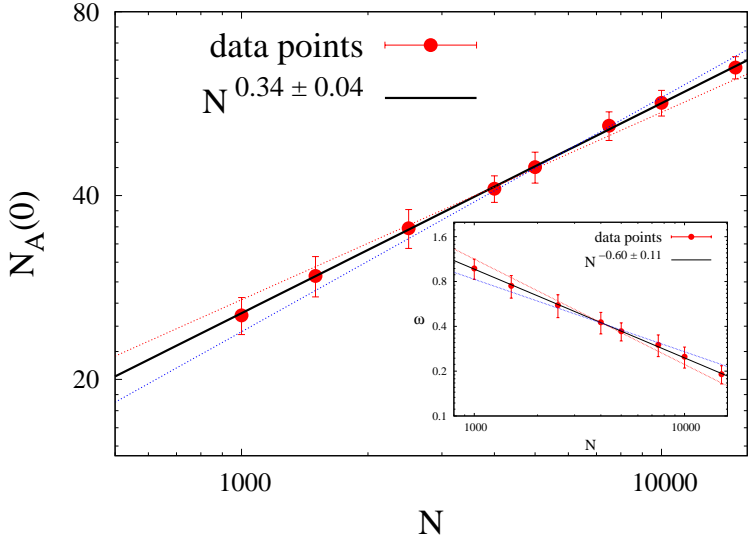}
	\caption{\label{fig:power_law} Plot shows that the $\langle N_A(0) \rangle$ scales as $N^{0.34 \pm 0.04}$, where $N$ is the total number of particles in the system. Inset of the figure shows that minimum width along the circumference of the circle which accommodates all the attractors at $t=0$ scales as $N^{-0.60 \pm 0.11}$}
\end{figure}

To quantify the results we carried out simulations for eight different densities in the range $\rho = 0.2-3.0$ (in units of $\sigma^{-2}$), keeping the radius $R$ of the circular disc fixed in all the simulation runs. For each density we ran the simulation for 50 different initial realizations (different random choises of particle positions within the circular disc) and computed the average number of attractors ($\langle N_A(0) \rangle$) at time $t=0$ in the system. Fig. \ref{fig:power_law} shows the dependence of $\langle N_A(0) \rangle$ on the total number of particles $N=\rho \pi R^2$ in the system (log-log plot). From the slope of the line we find that $\langle N_A(0) \rangle \propto N^{0.34\pm 0.04}$. Since the initial number of attractors is proportional to the area of the annular region of outer radius $R$ and width $\omega$ along the periphery where all the attractors lie ($\omega$ is averaged over all initial realizations), $\langle N_A(0)\rangle \propto \rho \omega$ or $\omega \propto \langle N_A(0)\rangle \rho^{-1}$. As $\rho\propto N$, we get $\omega \propto N^{-0.66 \pm 0.04}$.  Direct measurement of $\omega$ from the simulations yields $\omega \propto N^{-0.6 \pm 0.11}$(see inset of Fig. \ref{fig:power_law} ), consistent with the scaling obtained using $\langle N_A(0) \rangle $ simulation data. Thus, as the density of particles increases, the number of attractors increases as well but the width of the region in which these attractors lie decreases.
    
\begin{figure}[ht]
\includegraphics[width=0.6\linewidth]{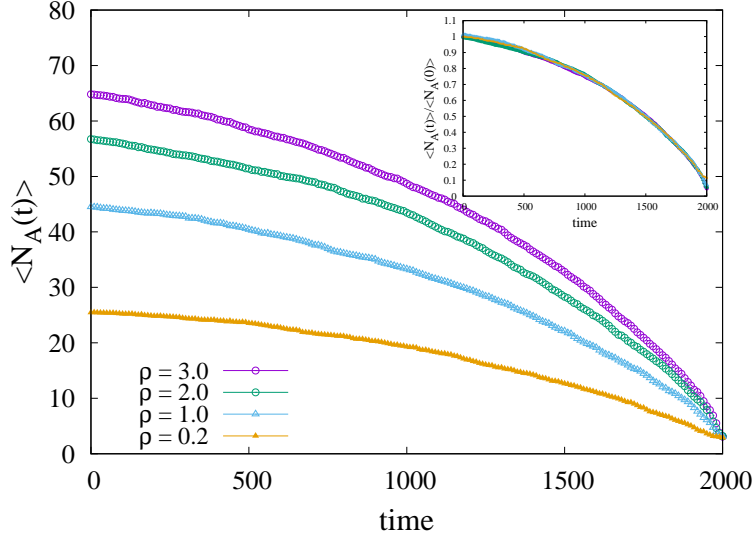}
	\caption{\label{fig:NTvst} $\langle N_T(t)\rangle$ vs $t$ plot for circular geometry of same radius and four different densities ($\rho = 0.2, 1.0, 2.0, 3.0$).  $\langle N_T(t)\rangle/\langle N_T(0)\rangle $ is plotted in the inset where data of all four densities collapses to give universal decay behavior of relative number of average attractors.}
\end{figure}

In Fig. \ref{fig:NTvst} we plot the average (over initial realizations) number of attractors, $\langle N_A(t)\rangle$ in the system as a function of time, $t$, for different particle densities. In the inset of Fig. \ref{fig:NTvst}, we collapse the different curves on a single universal plot for different values of the density, by normalizing $\langle N_A(t)\rangle$  by $\langle N_A(0)\rangle$. Therefore,    
\begin{eqnarray}
 \frac{\langle N_A(t)\rangle}{\langle N_A(0)\rangle} = f(t)
\label{eq:ft},
\end{eqnarray}
where $f(t)$ is a universal function of time that does not depend on the particle density.  Interestingly, $f(t)$ can not be fitted by a simple/stretched exponential or by power law decay and we have no analytic model for it.

\begin{figure}[h]
\includegraphics[width=\linewidth]{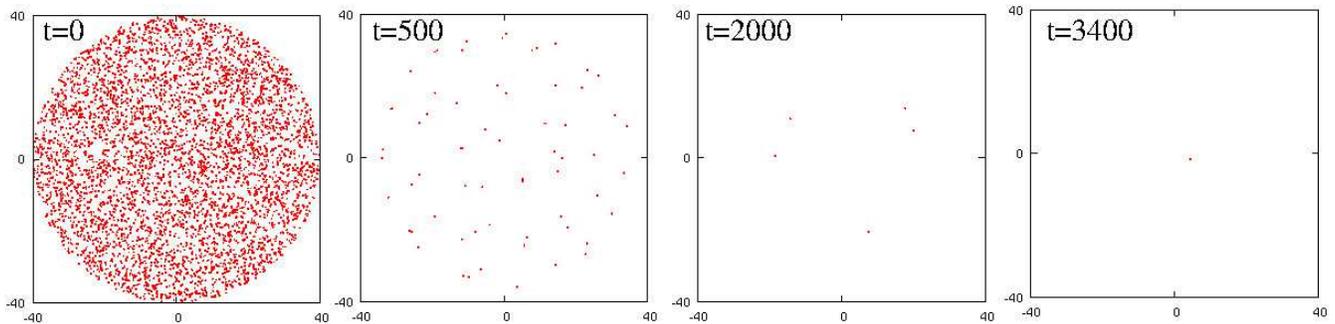}
	\caption{\label{fig:nearest_point} Snapshots of a system of 5000 particles at four different times where particles interact with their nearest neighbour. At $t=0$, particles are randomly distributed inside a circular region. Particles coalesce to form groups($t=500$), these groups again coalesce to form groups containing large number of particles($t=2000$) and process continues till all particles collapses to one group ($t=3400$).}
\end{figure}

The non-local dynamical rule we proposed can be used to assemble all the agents at a certain location in space and it is instructive to compare its efficiency to  that of a local algorithm according to which at every time step each particle finds the closest particle and moves towards it. In order to avoid short distance singularities, we introduce the constraint that if the separation between the particles is smaller than $\Delta x =0.02$, they stop sensing each other and each of them moves towards the next nearest particle located at distance larger than $0.02$. Repeated application of this algorithm results in the formation of numerous pointlike (of size $<\Delta x$) clusters of particles. These clusters coalesce to form new pointlike clusters composed of increasingly larger numbers of particles, a process reminiscent of nucleation and growth in phase separating systems \cite{Binder1977}. 
This process continue until all $N$ particles assemble into a single pointlike cluster of size $<\Delta x$. Fig. \ref{fig:nearest_point} shows the snapshots, at different times, of a system of $5000$ particles (initially randomly placed  in a circular disc with density $\rho=1/\sigma^2$) evolved using this algorithm (also see movie M9 in SI). We computed the total time of assembly of all the particles using this algorithm and compared it with the collapse time of our non-local algorithm. We found that the system assembles much faster using the non-local ($t=2200$)  than the local ($t=3400$) algorithm. Another point to notice is that the collapse zone is determined very early by the system using the proposed non-local algorithm and the particles always move towards the collapse zone. Particles which are closer to the collapse zone cover small distance, which progressively increases for the far particles, to reach the zone. Unlike this, the collapse zone is determined very late in the system using local algorithm and the motion of the particles is not always towards the collapse zone, therefore, all particles cover very large distance before reaching the collapse zone.     

In order to check whether and how the assembly of particles depends on dimensionality of the system, we performed simulations using the non-local algorithm in one and in three dimensions. In the 1D case, we randomly placed  $N=100$ particles on a line between $-L$ to $L$, where $2L=N/\rho_{1d}$ with uniform linear density $\rho_{1d}= 1/\sigma$. In this case the motion of all the particles is determined only by the two boundary particles that are the closest to ends of the line, i.e. to $-L$ and $L$. All particles which are on one side of the mid-point M of the line joining these two boundary particles, move along this line towards the particle closest to end of line on the other side of M. The two boundary particles remain the attractors of the system until the end of the collapse when all the particles reach the assembly zone around M (see Movie M10 in SI).  In 3D, we randomly placed $N=10000$ particles in a spherical region of radius $R=(3N/4\pi\rho_{3d})^{1/3}$ with uniform density $\rho_{3d}=1/\sigma^3$. Similar to the 2D case, we observed that the collapse proceeds through formation of lines that radiate outward from a point close to the center of the sphere (see Fig. S1 in the SI). Although visualization is more difficult in 3D than in 2D, we conclude that each of the lines is the common tangent to neighboring cones each of which contains the followers of a given attractor (not shown), that replace the slices shown in Fig.  \ref{fig:collapse1} for 2D case.

\section*{Conclusions}
In this paper, we simulated a ensemble of particles randomly distributed on a disc in two dimensions that follow a simple dynamical rule: every particle (follower) moves towards the farthest agent (attractor) from it. An obvious consequence of this dynamical rule is that the attractors are always located near the instantaneous outer boundary of the system and constitute a small fraction of the total number of particles $N$. As a follower moves towards its attractor, it approaches the perpendicular bisector of the imaginary line joining this attractor to its neighbouring attractor, and from this point on it executes a zigzag motion about this line as it switches between the two attractors; since the deviations from the line are small, it appears that the particle moves along the line. As time progresses, the system collapses but this collapse is anisotropic: the initially isotropic system self-organizes into slices of low particle density that are separated by lines of increasingly higher density and most particles move along these lines towards the assembly zone. We find that the initial number of attractors $\langle N_A(0)\rangle$ scales as $N^{0.34 \pm 0.4}$ and decreases with time as some of the attractors lose their followers and therefore forego their status of attractors; plotting the ratio $\langle N(t)\rangle$/$\langle N(0)\rangle$ vs $t$ yields a universal curve for all densities. We also found that line formation is not limited to circular disc geometry: lines are observed in square and semi-circular geometries as well, even though the number of lines in these geometries is much smaller and does not strongly depend on the initial conditions. Formation of lines in a circular disc geometry was observed for random non-uniformly distributed 2D particle systems as well, e.g. for radially non-uniform distribution (density varying as $1/r$) and hyper-uniform distribution\cite{Stillinger2003} and also in a non-random system in which particles were placed on a square lattice bounded by a circle. We found that collapse along lines is a unique feature of our non-local dynamical rule and appears to occur in 1, 2 and 3 dimensions. Even though the observation of such a collapse in 1D system appears to be trivial, it is actually not. For example, if one uses a local rule according to which particles move towards their nearest neighbors, the dynamics leads to the formation of many point-like clusters (each composed of several particles) which continue to coalesce until only one point-like cluster that contains all the particle in the system remains. This should be contrasted with the 1D dynamics produced by the non local rule where all the particles move uniformly towards the center of the 1D distribution. Interestingly, the time of assembly of a system evolved using the non-local rule is shorter than that for the local rule. While the non-local rule appears to be unphysical for most natural systems in which interactions are local in character, it can be implemented in artificial agent systems e.g., robots that can communicate across arbitrary distances. In recent years, self-organization and collective behavior in swarms of robots where a large number of robots interact using an simple algorithm has gained lot of attention \cite{Oprea,Sperati, Cao, Murakami, Trianni, Moeslinger, Baldassarre}.

\bibliography{references}

\section*{Acknowledgements}
We would like to thank Stas Burov, Yoav Soen, David Kessler, Baruch Barzel and Reuven Cohen for valuable discussions. This work was supported by grants from the Israel Science Foundation and from the Israeli Centers for Research Excellence program of the Planning and Budgeting Committee.

\section*{Author contributions}
K.S. conceived the idea and performed simulations. K.S. and Y.R. together analyzed the results and wrote the manuscript.

\section*{Additional information}
Competing Interests: The authors declare no competing interests.

\end{document}